\documentclass[10pt,shortpaper,twoside,web]{ieeeconf}
\usepackage{enumerate}
\usepackage{amssymb,epsfig,color,amsmath,amsfonts,mathrsfs,algorithm,algorithmic}
\usepackage{epstopdf}
\usepackage{mathrsfs,cite}
 \usepackage{graphicx} 
\usepackage{epstopdf}
\usepackage{subfigure}
\usepackage{booktabs}
\usepackage{graphicx}
\usepackage{stfloats}

\usepackage{latexsym,amsmath,amssymb,amsfonts,epsfig,graphicx,psfrag}
\usepackage{eepic,color,colordvi,amscd,mathtools,bm
}
\newtheorem{lemma}{Lemma}
\newtheorem{theorem}{Theorem}

\newtheorem{remark}{Remark}

\newtheorem{assumption}{Assumption}

\IEEEoverridecommandlockouts                              
\title{\LARGE \bf
Distributed Mirror Descent Algorithm with Bregman Damping for Nonsmooth Constrained Optimization}

\author{
	Guanpu Chen, Weijian Li, Gehui Xu, and Yiguang Hong, \IEEEmembership{Fellow, IEEE}
\thanks{This work was supported in part by Shanghai
	Municipal Science and Technology Major Project under Grant 2021SHZDZX0100, and in part by the National Natural Science Foundation
	of China under Grant 61733018. (\textit{Corresponding author: Yiguang Hong})}
\thanks{G. Chen and G. Xu is with Key Laboratory of Systems and Control, Academy of
	Mathematics and Systems Science, Beijing, China, and is also with School
	of Mathematical Sciences, University of Chinese Academy of Sciences,
	Beijing, China.
	{\tt\small chengp@amss.ac.cn, xghapple@amss.ac.cn}}%
\thanks{W. Li is with Department of Automation, University of Science and Technology of China, Hefei, Anhui, China.
        {\tt\small ustcwjli@mail.ustc.edu.cn}}
\thanks{Y. Hong is with Department of Control Science and Engineering \&
	Shanghai Research Institute for Intelligent Autonomous Systems, Tongji
	University, Shanghai, China, and is also with the Key Laboratory of Systems
	and Control, Academy of Mathematics and Systems Science, Beijing, China.
    	{\tt\small yghong@iss.ac.cn}}%
}

\begin{document}

\maketitle
\thispagestyle{empty}
\pagestyle{empty}

\begin{abstract}
To solve distributed optimization efficiently with various constraints and nonsmooth functions, we propose a distributed mirror descent algorithm with embedded Bregman damping, as a generalization of conventional distributed projection-based algorithms. In fact, our continuous-time algorithm well inherits good capabilities of mirror descent approaches to rapidly compute explicit solutions to the problems with some specific constraint structures. Moreover, we rigorously  prove the convergence of our algorithm, along with the boundedness of the trajectory and the accuracy of the solution.


\end{abstract}

\section{INTRODUCTION}
Distributed optimization has served as a hot topic in recent years for its broad applications in various fields such as sensor networks and smart grids \cite{yang2019survey,zhu2011distributed,yuan2015regularized,cherukuri2016initialization,xu2016distributed,you2018distributed,lu2019online}. 
Under multi-agent frameworks, the global cost function
consists of agents' local ones, and each
agent shares limited information with its neighbors
through a network to achieve an optimal solution. Meanwhile, distributed continuous-time algorithms have been well developed thanks to system dynamics and control theory \cite{yang2016multi,zhu2018continuous,liang2017distributed,li2019distributed,li2021exponentially}.

Up to now, various approaches have been employed for distributed design of  constrained optimization. Intuitively, implementing projection operations on local constraints is a most popular method, such as projected proportional-integral protocol \cite{yang2016multi} and projected dynamics with constraints based on KKT conditions \cite{zhu2018continuous}. In addition, other approaches such as  distributed primal-dual dynamics and penalty-based algorithms \cite{liang2017distributed,li2019distributed,li2021exponentially} also perform well provided that constraints are endowed with specific expressions.
However, time complexity in finding optimal solutions with complex or high-dimensional constraints forces researchers to exploit efficient approaches for special constraint structures, such as the unit simplex and the Euclidean sphere. 

In fact, the mirror descent (MD) method serves as a powerful tool for solving constrained optimization. As we know, first introduced in \cite{nemirovskij1983problem}, MD is regarded as a generalization of (sub)gradient methods. With mapping the variables into a conjugate space, MD employs the Bregman divergence and performs well in handling local constraints with specific structures \cite{krichene2015accelerated,diakonikolas2019approximate,ben2001ordered}. This process results in a faster convergence rate than that of projected (sub)gradient descent algorithms, especially for large-scale optimization problems. Undoubtedly, as such an important tool, MD has played a crucial role in various distributed algorithm designs, as given in \cite{shahrampour2017distributed,yuan2018optimal,sun2020distributed,wang2020distributed}.
  
In recent years, continuous-time MD-based algorithms have also attracted much attention. For example,  \cite{krichene2015accelerated}  proposed the acceleration of a continuous-time MD algorithm, and afterward, \cite{xu2018continuous} showed continuous-time stochastic MD for strongly convex functions, while \cite{gao2020continuous} proposed a discounted continuous-time MD dynamics to approximate the exact solution. In the distributed design, 
although \cite{sun2020distributed} presented a distributed MD dynamics with integral feedback, the result merely achieved optimal consensus and part variables turn to be unbounded. Due to the booming development and extensive demand of distributed design, distributed continuous-time MD-based methods need more exploration and development actually.

Therefore, we study continuous-time MD-based algorithms to solve distributed nonsmooth optimization with local and coupled constraints. 
The main contributions of this note can be summarized as follows. We propose a distributed continuous-time mirror descent algorithm by introducing the Bregman damping, which can be regarded as a generalization of classic distributed projection-based dynamics \cite{yang2016multi,zeng2018distributed} by taking the {Bregman damping} in a quadratic form. 
Moreover, our algorithm well inherits the good capabilities of MD-based approaches to rapidly compute explicit solutions to the problems with some concrete constraint structures like the unit simplex or the Euclidean sphere. With the designed {Bregman damping}, our MD-based algorithm makes all the variables' trajectories bounded, which could not be ensured in \cite{krichene2015accelerated,sun2020distributed}, and avoids the inaccuracy of the convergent point occurred in \cite{gao2020continuous}.

The remaining part is organized as follows. Section \ref{Sec:Pre} gives related preliminary knowledge. Next, Section \ref{Sec:for} formulates the distributed optimization and provides our algorithm, while Section \ref{Sec:main} presents the main results. Then, Section \ref{Sec:num} provides illustrative numerical examples. Finally, Section \ref{Sec:con} gives the conclusion.
\section{Preliminaries}\label{Sec:Pre}
In this section, we give necessary notations and related preliminary knowledge.
\subsection{Notations}
Denote $\mathbb R^n$ (or $\mathbb R^{m\times n}$) as the set of $n$-dimensional (or $m$-by-$n$) real column vectors (or real matrices),  and $I_n$ as the $n\times n$ identity matrix. 
Let $ {1}_n$ (or ${0}_{n}$) be the $n$-dimensional column vector with all entries of $1$ (or $0$). Denote $A \otimes B$ as the Kronecker product of matrices $A$ and $B$. Take $col\{x_1,\dots,x_n\}=col\{x_i\}_{i=1}^n=(x_1^{\rm T},\dots,x_n^{\rm T})^{\rm T}$, $\|\cdot\|$  as the Euclidean norm, and $\text{rint}(C)$ as the relative interior of the set $C$ \cite{boyd2004convex}. 

An undirected graph can be defined by $\mathcal G(\mathcal V, \mathcal E)$, where $\mathcal V=\{1,\ldots, n\}$ is the set of nodes and $\mathcal E\subset\mathcal V \times \mathcal V$ is the set of edges. Let $ \mathcal A=[a_{ij}]\in\mathbb R^{n\times n}$ be the {\em adjacency matrix} of $\mathcal G$ such that $a_{ij}=a_{ji}>0$ if $\{j,i\}\in\mathcal E$, and $a_{ij}=0$, otherwise. The {\em Laplacian matrix} is $L_n=\mathcal D- \mathcal A$, where $\mathcal D=\text{diag}\{D_{ii}\}\in\mathbb R^{n\times n}$ with $\mathcal D_{ii}=\sum_{j=1}^n a_{ij}$.
If the graph $\mathcal G$ is connected, then $\text{ker}(L_n)=\{k{1}_n:k\in\mathbb R\}$.
\subsection{Convex analysis}
For a closed convex set $\Omega\subseteq\mathbb R^n$, the projection map $P_\Omega:\mathbb R^n\rightarrow\Omega$ is defined as $P_\Omega(x)=\text{argmin}_{y\in\Omega}\|x-y\|$. Especially, denote $[x]^+\triangleq P_{\mathbb R^n_+}(x)$ for convenience.
For $x\in\Omega$,  denote the normal cone to $\Omega$ at $x$ by
$$\mathcal N_\Omega(x)=\big\{v\in\mathbb R^{n}: v^{\rm T}(y-x)\leq 0,\quad\forall y\in\Omega\big\}.$$
A continuous function $f:\mathbb R^{n}\rightarrow \mathbb R$ is $\omega$-strongly convex on $\Omega$ if
$$(x-y)^{\rm T}(g_x-g_y)\geq\omega\|x-y\|^2,\,\forall x,y \in C,$$ 
where $\omega>0$, $g_x\in\partial f(x)$, and $g_y\in\partial f(y)$.

The Bregman divergence based on the differentiable generating function $f:\Omega\rightarrow \mathbb R$  is defined as
\begin{align*}
	D_f(x,y)=f(x)-f(y)-\nabla f(y)^T(x-y), \, \forall x,y\in \Omega.
\end{align*}
The convex conjugate function of $f$ is defined as
\begin{align*}
f^*(z)=\sup_{x\in \Omega}\big\{x^Tz-f(x)\big\}.
\end{align*}
The following lemma reveals a classical conclusion about convex conjugate functions, of which readers can find more details in  \cite{nemirovskij1983problem,diakonikolas2019approximate}.
\begin{lemma}\label{lemma2}
Suppose that a function $f$ is differentiable and strongly convex on a closed convex set $\Omega$. Then  $f^*(z)$ is convex and differentiable, and $f^*(z)=\min_{x\in \Omega}\big\{-x^Tz+f(x)\big\}$.  Moreover, $\nabla f^*(z)=\Pi_{\Omega}^f(z)$, 
where
	\begin{align}\label{Pi}
\Pi_{\Omega}^f(z)\triangleq\mathop{\text{argmin}}_{x\in \Omega}\big\{-x^Tz+f(x)\big\}.
\end{align} 
\end{lemma}
\subsection{Differential inclusion}
A differential inclusion is given by
\begin{equation}\label{dif_in}
\dot x(t)= \mathcal F(x(t)),
~x(0)=x_0,~t\ge 0,
\end{equation}
where $\mathcal F:\mathbb R^n \rightrightarrows \mathbb R^n$ is a set-valued map. $\mathcal F$ is upper semi-continuous at $x$ if  there exists $\delta>0$ for all $\epsilon>0$ such that 
\begin{equation*}
\mathcal F(y) \subset \mathcal F(x)+B(0;\epsilon), ~\forall y\in B(x;\delta)
\end{equation*}
and it is upper semi-continuous if it is so for all $x \in \mathbb R^n$.
A Caratheodory solution to (\ref{dif_in}) defined on $[0,\tau) \subset [0,+\infty)$ is an absolutely continuous function $x:[0,\tau) \rightarrow \mathbb R^n$ satisfying (\ref{dif_in}) for almost all $t\in [0,\tau)$ in Lebesgue measure \cite{cortes2008discontinuous}. The solution $x(t)$ is right maximal if it has no extension in time. 
A set $\mathcal M$ is said to be weakly (strongly) invariant with respect to (\ref{dif_in}) if $\mathcal M$ contains a (all) maximal solution to (\ref{dif_in}) for any $x_0\in \mathcal M$.
If $ 0_n \in \mathcal F(x_e)$, then $x_e$ is an equilibrium point of (\ref{dif_in}).
The existence of a solution to (\ref{dif_in}) is guaranteed by the following lemma \cite{cortes2008discontinuous}.
\begin{lemma}\label{lem:exist}
If $\mathcal F$ is locally bounded, upper semicontinuous, and takes nonempty, compact and convex values, then there exists
a Caratheodory solution to (\ref{dif_in}) for any initial value.
\end{lemma}

Let $V:\mathbb R^n \rightarrow \mathbb R$ be a locally Lipschitz continuous function, and $\partial V(x)$ be the 
Clarke generalized gradient of $V$ at $x$. The set-valued Lie derivative for $V$ is defined by
$\mathcal L_{\mathcal F} V(x) \triangleq\{a\in \mathbb R: a = p^T v, p\in \partial V(x),v\in \mathcal F(x)\}$. Let $\max \mathcal L_{\mathcal F} V(x)$ be the largest element of  $\mathcal L_{\mathcal F} V(x)$.
Referring to \cite{cortes2008discontinuous}, we have the following invariance principle for (\ref{dif_in}).
\begin{lemma}\label{lem:inv_prin} 
Suppose that $\mathcal F$ is upper semi-continuous and locally bounded, and $\mathcal F(x)$ takes nonempty, compact, and convex values. Let $V:\mathbb R^n \rightarrow \mathbb R$ be a locally Lipschitz and regular function, $\mathcal S \subset \mathbb R^n$ be compact and strongly invariant for (\ref{dif_in}), and $\psi(t)$ be a solution to (\ref{dif_in}). Take
$$\mathcal R=\{x\in \mathbb R^n:
0 \in \mathcal L_{\mathcal F}V(x)\},$$
and $\mathcal M$ as the largest weakly invariant subset of $\bar{\mathcal R} \cap \mathcal S$, where $\bar{\mathcal R}$ is the closure of $\mathcal R$. If $\max \mathcal L_{\mathcal F}V(x) \le 0$ for all $x\in \mathcal S$, then 
$\lim_{t\to \infty}{\rm dist}(\psi(t), \mathcal M)=0$.
\end{lemma}

\section{Formulation and algorithm}\label{Sec:for}

In this paper, we consider a nonsmooth optimization problem with both local and coupled constraints. There are $N$ agents indexed by $\mathcal V=\{1,\dots,N\}$ in a network $\mathcal G(\mathcal V,\mathcal E)$. For agent $i$, the decision variable is $x_i$, the local feasible set is $\Omega_i\subseteq \mathbb R^n$, and the local cost function is $f_i:\mathbb R^n\rightarrow \mathbb R$. Define $\bm\Omega=\prod_{i=1}^{N}\Omega_i$ and $\bm x = \text{col}\{x_i\}_{i=1}^N$.  All agents cooperate to solve the following distributed optimization:
\begin{align}\label{formulation}
\min_{\bm x\in\bm\Omega}&\;\sum_{i=1}^{N}f_i(x_i)\nonumber\\
\text{s.t.}&\;\sum_{i=1}^N g_i(x_i)\leq 0_p,\quad \sum_{i=1}^N A_ix_i- b_i=0_q,
\end{align}
where $g_i:\mathbb R^{n}\rightarrow \mathbb R^p$, $A_i\in\mathbb R^{q\times n}$  and $ b_i\in\mathbb R^q$, for $i\in\mathcal V$. Except for the local constraint $\bm\Omega$, other constraints in \eqref{formulation} are said to be coupled ones since the solutions rely on global information.
In the multi-agent network, agent $i$ only has the local decision variable $x_i\in\Omega_i$, and moreover, the local information $f_i$, $g_i$, $A_i$ and $b_i$. Thus, agents need communication with neighbors through the network $\mathcal G$.

Actually, MD replaces the Euclidean regularization in (sub)gradient descent algorithms with Bregman divergence. In return, different generating functions of Bregman divergence may efficiently bring explicit solutions on different special feasible sets.
For example, if $\phi(x)=\frac{1}{2}\|x\|^2$ and $\Omega$ is convex and closed, then 
\begin{align}\label{quadratic}
\Pi_{\Omega}^\phi(z)=
\mathop{\text{argmin}}_{x\in \Omega}\frac{1}{2}\|x-z\|^2=P_\Omega(z),
\end{align}
which actually turns into the classical Euclidean regularization with projection operations. Furthermore, if  $\Omega=\{x\in\mathbb R^n_+:\sum_{k=1}^nx^k=1\}$ and  $\phi(x)=\sum_{k=1}^nx^k\log(x^k)$ with the convention $0\log0=0$, then
\begin{align}\label{qq}
\Pi_{\Omega}^\phi(z)=\text{col}\Big\{\frac{\exp(z^k)}{\sum_{j=1}^n\exp(z^j)}\Big\}_{k=1}^n,
\end{align}
which is the well-known  KL-divergence on the unit simplex. 


Assign a generating function $\phi_i:\mathbb R^n\rightarrow \mathbb R$ of the Bregman divergence to each agent $i\in\mathcal V$. 
Then we consider the following assumptions for \eqref{formulation}.
\begin{assumption}\label{ass_1}
	\
	\begin{enumerate}[(i)]
		\item For $i\in\mathcal V$, $\Omega_i$ is closed and convex,  $ f_i$ and $g_i$ are convex on $\Omega_i$, and moreover,  $\phi_i$ is differentiable and strongly  convex on  $\Omega_i$.
		\item There exists at least one $\bm x\in\text{rint}(\bm \Omega)$ such that $\sum_{i=1}^{N}g_i(x_i)<  0_p$ and $\sum_{i=1}^{N}A_ix_i- b_i=0_q$.
		\item The undirected graph $\mathcal G$ is connected.
	\end{enumerate}
\end{assumption}

\begin{remark}
	Clearly, \eqref{formulation}  can be regarded as a generalization for both distributed optimal consensus problems \cite{shi2012randomized,li2021exponentially,sun2020distributed} and distributed resource  allocation problems \cite{yi2016initialization,zeng2018distributed}. Moreover,  $g_i$ in the coupled constraints may not required to be affine, which is more general than the constraints in previous works \cite{liang2017distributed,chen2021distributed}. Also, the problem setting does not require strongly or strictly convexity for either cost functions $f_i$ or constraint functions $g_i$ \cite{yi2016initialization,liang2017distributed}, and the selection qualification for generating function $\phi_i$ has also been widely used \cite{sun2020distributed,krichene2015accelerated,gao2020continuous}. 
\end{remark}

For designing a distributed algorithm, we introduce auxiliary variables $\omega_i\in\mathbb R^p$, $\nu_i\in\mathbb R^q$, $\lambda_i\in\mathbb R^p$, $\mu_i\in\mathbb R^q$, $y_i\in\mathbb R^n$ and $\gamma_i\in\mathbb R^p$ for each agent $ i\in\mathcal V$. Moreover, we employ the gradient $\nabla\phi_i$ of generating functions as the \textit{Bregman damping} in the algorithm, which ensures the trajectories' boundedness \cite{krichene2015accelerated,sun2020distributed} and avoids the convergence inaccuracy \cite{gao2020continuous}. 
Recall that $a_{ij}$ is the $(i,j)$-th entry of the adjacency matrix $\mathcal A$ and $\Pi_{\Omega_i}^{\phi_i}(\cdot)$ is defined in \eqref{Pi}. Then we propose a distributed  mirror descent algorithm with Bregman damping (MDBD) for \eqref{formulation}.

\begin{algorithm}
	\caption{MDBD for $i\in\mathcal V$}
	\label{alg:1}
	\quad\\[1pt]
	\textbf{Initialization:}
	\\[2pt]
	$x_i(0)\in\Omega_i$, $y_i(0)=0_n$, $\lambda_i(0)=0_p$, $\gamma_i(0)=0_p$, 
	\\$\omega_i(0)=0_p$, $\mu_i(0)=0_q$, $\nu_i(0)=0_q$;
	\\take a proper generating function $\phi_i(\cdot)$ according to $\Omega_i$.
	\\[2pt]
	\textbf{Flows renewal:}
	\begin{flalign*}
		\dot y_i\,&{\in}- { \partial}f_i(x_i)-{ \partial} g_i(x_i)^{\rm T}\lambda_i -A_i^{\rm T}\mu_i+\nabla \phi_i(x_i)-y_i,\\
\dot\gamma_i&= g_i(x_i)-\sum_{j=1}^Na_{ij}(\omega_i-\omega_j)+\lambda_i-\gamma_i,\\
\dot\mu_i &= A_ix_i-b_i-\sum_{j=1}^Na_{ij}(\nu_i-\nu_j),\\
	\dot\omega_i & = \sum_{j=1}^Na_{ij}(\lambda_i-\lambda_j),\\
	\dot\nu_i & = \sum_{j=1}^Na_{ij}(\mu_i-\mu_j),\\
	x_i&= \Pi_{\Omega_i}^{\phi_i}(y_i),\\
	\lambda_i&=[\gamma_i]^+.
	\end{flalign*}
\end{algorithm}
In Algorithm \ref{alg:1}, for each agent $i\in\mathcal V$, information like $ { \partial}f_i$, $ { \partial} g_i$, $A_i$ and $b_i$ serves as private knowledge, and values like $\omega_i$, $\nu_i$, $\lambda_i$ and  $\mu_i$ should be exchanged with neighbors through the network $\mathcal G$. Moreover, generating function $\phi_i$ and {Bregman damping $\nabla\phi_i$} can be determined privately and individually, not necessarily identical.  It follows from Lemma \ref{lem:exist} that
the existence of a Caratheodory solution to Algorithm \ref{alg:1} can be guaranteed.

For simplicity, define
$$\bm \lambda=\text{col}\big\{\lambda_i\big\}_{i=1}^N\in\mathbb R^{Np},~~\bm \mu=\text{col}\big\{\mu_i\big\}_{i=1}^N\in\mathbb R^{Nq},$$
$$\bm \omega=\text{col}\big\{\omega_i\big\}_{i=1}^N\in\mathbb R^{Np},~~\bm \nu=\text{col}\big\{\nu_i\big\}_{i=1}^N\in\mathbb R^{Nq},$$
$$\bm y=\text{col}\big\{y_i\big\}_{i=1}^N\in\mathbb R^{Nn},~~\bm \gamma=\text{col}\big\{\gamma_i\big\}_{i=1}^N\in\mathbb R^{Np}.$$
Let  $\bm \Theta=\bm\Omega\times\mathbb R^{Np}_+\times \mathbb R^{N(p+2q)}$, and moreover, 
$$\bm z=\text{col}\{\bm x,\bm \lambda,\bm \mu,\bm \omega,\bm \nu\},\quad\bm s=\text{col}\{\bm y,\bm \gamma,\bm \mu,\bm \omega,\bm \nu\}.$$
Take the Lagrangian function $\mathcal L:\bm\Theta\rightarrow \mathbb R$ as
\begin{align}\label{lagrangian}
\mathcal L(\bm z)=&\sum_{i=1}^Nf_i(x_i)+\sum_{i=1}^N\lambda_i^T\big(g_i(x_i)-\sum_{j=1}^Na_{ij}(\omega_i-\omega_j)\big)\nonumber\\
&+\sum_{i=1}^N\mu_i^T\big(A_ix_i-b_i-\sum_{j=1}^Na_{ij}(\nu_i-\nu_j)\big).
\end{align}
For such a distributed convex optimization \eqref{formulation} with a zero dual gap, $\bm x^\star\in\bm\Omega$ is an optimal solution to problem \eqref{formulation} if and only if there exist auxiliary variables $(\bm \lambda^\star,\bm \mu^\star,\bm \omega^\star,\bm \nu^\star)\in\mathbb R^{Np}_+\times\mathbb R^{N(p+2q)}$, such that $\bm z^\star=\text{col}\{\bm x^\star,\bm \lambda^\star,\bm \mu^\star,\bm \omega^\star,\bm \nu^\star\}$ is a saddle point of $\mathcal L$ \cite{zeng2018distributed}, that is,
for arbitrary $\bm z\in\bm \Theta$,
\begin{align*}
\mathcal L(\bm x^\star,\bm \lambda,\bm \mu,\bm \omega^\star,\bm \nu^\star)\leq\mathcal L(\bm x^\star,\bm \lambda^\star,&\bm \mu^\star,\bm \omega^\star,\bm \nu^\star)\nonumber\\
\leq&\mathcal L(\bm x,\bm \lambda^\star,\bm \mu^\star,\bm \omega,\bm \nu).
\end{align*}
Define
\begin{align}\label{gradient}
F(\bm z)=
\begin{bmatrix}
\text{col}\big\{ {\partial} f_i(x_i)+ { \partial} g_i(x_i)^{\rm T}\lambda_i +A_i^{\rm T}\mu_i\big\}_{i=1}^N\\
\text{col}\big\{-g_i(x_i)\big\}_{i=1}^N+ \bm L_p\bm w\\
\text{col}\big\{-A_ix_i+b_i\big\}_{i=1}^N+ \bm L_q\bm \nu\\
-\bm L_p\bm \lambda\\
-\bm L_q\bm \mu
\end{bmatrix},
\end{align}
where $\bm L_{p}=L_N\otimes I_{p}$ and $\bm L_q=L_N\otimes I_{q}$.
In fact, $\bm z^\star$ is a saddle point of $\mathcal L$ if and only if $-F(\bm z^\star)\in\mathcal N_{\bm \Theta}({\bm z}^\star)$, which was obtained in \cite{yi2016initialization,zeng2018distributed,liang2017distributed}.

Hence, Algorithm \ref{alg:1} can be presented in the following compact form:
\begin{equation}\label{compact}
\left\{
\begin{aligned}
\dot {\bm s}\;{\in}&-F(\bm z)+\nabla\Phi(\bm z)-\bm s,\\
\bm z=&\; \Pi_{\bm\Theta}^\Phi(\bm s),
\end{aligned}
\right.
\end{equation}
where 
\begin{subequations}
\begin{align}
\nabla\Phi(\bm z)\triangleq&\text{col}\Big\{\text{col}\{\nabla\phi_i(x_i)\}_{i=1}^N,\text{col}\{\lambda_i\}_{i=1}^{N},\nonumber\\
&\quad \text{col}\{\mu_i\}_{i=1}^{N},\text{col}\{\omega_i\}_{i=1}^{N},\text{col}\{\nu_i\}_{i=1}^{N}\Big\},\label{notation1}\\
\Pi_{\bm\Theta}^\Phi(\bm s)\triangleq&\text{col}\Big\{\text{col}\{\Pi_{\Omega_i}^{\phi_i}(y_i)\}_{i=1}^N,\text{col}\{[\gamma_i]^+\}_{i=1}^{N},\nonumber\\
&\quad \text{col}\{\mu_i\}_{i=1}^{N},\text{col}\{\omega_i\}_{i=1}^{N},\text{col}\{\nu_i\}_{i=1}^{N}\}\Big\}.\label{notation2}
\end{align} 
\end{subequations}

\begin{remark}
In fact, if $\phi(\cdot)=\frac{1}{2}\|\cdot\|^2$, then $\Pi_{\mathbb R^p}^\phi(z)=P_{\mathbb R^p}(z)=z$ and $\Pi_{\mathbb R^p_+}^\phi(z)=P_{\mathbb R^p_+}(z)=[z]^+$, and therefore, \eqref{compact} can be rewritten as 
	\begin{equation}\label{pro_alg}
	\left\{
	\begin{aligned}
	\dot {\bm s}\;{ \in}&-F(\bm z)+\bm z-\bm s,\\
	\bm z=&\; P_{\bm\Theta}(\bm s),
	\end{aligned}
	\right.
	\end{equation}
	which is actually a widely-investigated dynamics such as the proportional-integral protocol in \cite{yang2016multi} and projected output feedback in \cite{zeng2018distributed}.   
	Thus, MDBD generalizes the conventional distributed projection-based design for constrained optimization. Obviously, $\bm z$ in \eqref{pro_alg} is replaced with the {Bregman damping $\nabla\Phi(\bm z)$} in \eqref{compact}.
\end{remark}

\section{Main results}\label{Sec:main}
In this section, we investigate the convergence of MDBD.  Though Bregman damping improves the convergence of MDBD, the process also brings challenges for the convergence analysis. The following lemma shows the relationship between MDBD and the saddle points of the Lagrangian function  $\mathcal L$.
\begin{lemma}\label{equivelance}
Under Assumption \ref{ass_1}, $\bm z^\star$ is a saddle point of Lagrangian function $\mathcal L$ in \eqref{lagrangian} if and only if there exists $\bm s^\star\in-F(\bm z^\star)+\nabla \Phi(\bm z^\star)$ such that $\bm z^\star=\Pi_{\bm \Theta}^{\Phi}(\bm s^\star)$.
\end{lemma}
\textbf{Proof}.
For $\tilde {\bm z}=\Pi_{\bm \Theta}^{\Phi}(\bm s^\star)$,  the first-order condition is 
\begin{align}\label{firstorder_Pi}
-F(\bm z^\star)+\nabla\Phi(\bm z^\star)-\nabla\Phi(\tilde{\bm z})\in\mathcal N_{\bm \Theta}(\tilde{\bm z}).
\end{align}
We firstly show the sufficiency. Given $\bm z^\star$, suppose that  there exists $\bm s^\star\in-F(\bm z^\star)+\nabla \Phi(\bm z^\star)$ such that $\bm z^\star=\Pi_{\bm \Theta}^{\Phi}(\bm s^\star)$.
Thus, \eqref{firstorder_Pi} holds with $\tilde{\bm z}=\bm z^\star$, and $-F(\bm z^\star)\in\mathcal N_{\bm \Theta}({\bm z}^\star)$, which means that $\bm z^\star$ is  a saddle point of  $\mathcal L$.

Secondly, we show the necessity. Suppose  $-F(\bm z^\star)\in\mathcal N_{\bm \Theta}({\bm z}^\star)$ and take $\bm s^\star\in-F(\bm z^\star)+\nabla \Phi(\bm z^\star)$. Recall that \eqref{firstorder_Pi} holds with $\tilde{\bm z}=\bm z^\star$, which implies that $\bm z^\star$ is a solution to $\Pi_{\bm\Theta}^\Phi(\bm s^\star)$. Furthermore, since $\phi_i(\cdot)$ and $\frac{1}{2}\|\cdot\|^2$ are strongly convex, the solution to $\Pi_{\bm\Theta}^\Phi(\bm s^\star)$ is unique. Therefore, $\bm z^\star=\Pi_{\bm\Theta}^\Phi(\bm s^\star)$. \hfill$\square$

The following theorem shows the correctness and the convergence of  Algorithm \ref{alg:1}.
\begin{theorem}
	 Under Assumption \ref{ass_1}, the following statements hold.
	\begin{enumerate}[(i)]
		\item The trajectory $(\bm s(t),\bm z(t))$ of \eqref{compact} is bounded;
		\item $\bm x(t)$ converges to an optimal solution to problem \eqref{formulation}.
	\end{enumerate}
\end{theorem}
\textbf{Proof}.
(i) Firstly, we show that the output $\bm z(t)$ is bounded. By Lemma \ref{equivelance}, take $\bm z^\star$ as a saddle point of $\mathcal L$ and thus, there exists $\bm s^\star\in-F(\bm z^\star)+\nabla \Phi(\bm z^\star)$ such that $\bm z^\star=\Pi_{\bm \Theta}^{\Phi}(\bm s^\star)$. Take $\phi^*_i$ as the convex conjugate of $\phi_i$, and construct a Lyapunov candidate function as 
\begin{align}\label{Lya_1}
V_1=&\sum_{i=1}^ND_{\phi_i^*}( y_i-y_i^\star)+\frac{1}{2}\|\bm \gamma-\bm \gamma^\star\|^2+\frac{1}{2}\|\bm \mu-\bm \mu^\star\|^2\nonumber\\
&+\frac{1}{2}\|\bm \omega-\bm \omega^\star\|^2+\frac{1}{2}\|\bm \nu-\bm \nu^\star\|^2.
\end{align}
Since $x_i= \Pi_{\Omega_i}^{\phi_i}(y_i)$, it follows from Lemma \ref{lemma2} that
\begin{subequations}\label{fenchel}
\begin{align}
\phi^*_i(y_i)=&x_i^Ty_i-\phi_i(x_i),\\
\phi^*_i(y_i^\star)=&x_i^{\star T}y_i^\star-\phi_i(x_i^\star).
\end{align}
\end{subequations}
Thus, by substituting \eqref{fenchel}, the Bregman divergence becomes
\begin{align*}
D_{\phi_i^*}(y_i-y_i^\star)=&\phi^*_i(y_i)-\phi^*_i(y_i^\star)-\nabla \phi^*_i(y_i^\star)^T(y_i-y_i^\star)\\
=&\phi_i(x_i^\star)-\phi_i(x_i)-(x_i^{\star }-x_i)^Ty_i.
\end{align*}
Since $\phi_i(\cdot)$ is strongly convex for $i\in\mathcal V$, there exists a positive constant $\sigma$ such that
\begin{align*}
\sum_{i=1}^N\!D_{\phi_i^*}\!( y_i\!-\!y_i^\star)\!\geq\! \frac{\sigma}{2}\!\|\bm x\!-\!\bm x^\star\!\|^2\!+\!\sum_{i=1}^N( x_i^{\star }\!- \!x_i)^T\!(\nabla \phi_i( x_i)\!-\!y_i).
\end{align*}
In fact, $\nabla\phi^*_i(y_i)=\mathop{\text{argmin}}_{x\in {\Omega_i}}\{-x^Ty_i+\phi_i(x)\}$, which leads to
\begin{align*}
0\leq&\big(\nabla\phi_i(\nabla\phi^*_i(y_i))-y_i\big)^T\big(\nabla \phi^*_i(y_i^\star)-\nabla\phi_i^*(y_i)\big)\\
=&(\nabla \phi_i(x_i)-y_i)^T(x_i^{\star }-x_i).
\end{align*}
Thus, 
$\sum_{i=1}^ND_{\phi_i^*}( y_i-y_i^\star)\geq \frac{\sigma}{2}\|\bm x-\bm x^\star\|^2$. In addition,
\begin{align*}
\|\bm\lambda-\bm\lambda^\star\|^2=\|[\bm \gamma]^+-[\bm\gamma^\star]^+\|^2\leq\|\bm \gamma-\bm\gamma^\star\|^2.
\end{align*}
Therefore, 
\begin{align}\label{V_1radially_unbdd}
V_1(\bm s(t))\geq&\frac{\kappa}{2}\Big(\|\bm x-\bm x^\star\|^2+\|\bm \lambda-\bm \lambda^\star\|^2+\|\bm \mu-\bm \mu^\star\|^2\nonumber\\
&+\|\bm \omega-\bm \omega^\star\|^2+\|\bm \nu-\bm \nu^\star\|^2\Big),
\end{align}
where $\kappa=\min\{\sigma,1\}$. This means $V_1(\bm s(t))\geq\frac{\kappa}{2}\|\bm z-\bm z^\star\|^2$, that is, $V_1$ is radially unbounded  in $\bm z$. 
Clearly, the function $V_1$ along (20) satisfies
\begin{equation}	
\begin{aligned}\label{deriv_V1}
\mathcal L_{\mathcal F}&V_1=\big\{\beta\in \mathbb R, \beta=\sum_{i=1}^N\big(\nabla\phi_i^*(y_i)-\nabla\phi_i^*(y_i^\star)\big)^T\dot y_i \\
+&(\bm \gamma \!-\!\bm \gamma^\star)^T\dot{\bm \gamma}\nonumber
+(\bm \mu \!-\!\bm \mu^\star)^T\dot{\bm \mu}
+(\bm \omega \!-\!\bm \omega^\star)^T\dot{\bm \omega}+(\bm \nu\!-\!\bm \nu^\star)^T\dot{\bm \nu} \\
=&\big\{\beta\in \mathbb R, \beta=\big(\bm z\!-\!\bm z^\star\big)^T\!\big(-\bm \eta\!+\!\nabla \Phi(\bm z)\!-\!\bm s\big),\; \bm \eta\in F(\bm z)\big\}.
\end{aligned}
\end{equation}
Combining the convexity of $f_i$ and $g_i$ with the property for saddle point $\bm z^\star$,
\begin{align}\label{optimal}
-(\bm z-&\bm z^\star\big)^T\bm \eta\\
\leq& \mathcal L(\bm x^\star,\bm \lambda,\bm \mu,\bm \omega^\star,\bm \nu^\star)-\mathcal L(\bm x,\bm \lambda^\star,\bm \mu^\star,\bm \omega,\bm \nu)\leq 0.\nonumber
\end{align}
Thus, 
\begin{align}\label{deriv_V11}
{\beta}&\leq (\bm z-\bm z^\star)^T(\nabla \Phi(\bm z)-\bm s)\\
&=\sum_{i=1}^N(x_i-x_i^\star)^T(\nabla\phi_i(x_i)-y_i)+(\bm \lambda-\bm \lambda^\star)^T(\bm \lambda-\bm \gamma).\nonumber
\end{align}
On the one hand, for $i\in\mathcal  P$, we consider a differentiable function 
\begin{align*}
J(\alpha)=\phi_i\big(\alpha x_i^\star+(1-\alpha)x_i\big)-\big(\alpha x_i^\star+(1-\alpha)x_i\big)^Ty_i,
\end{align*}
with a constant $\alpha\in[0,1]$. Correspondingly, we have
\begin{align*}
J'(\alpha)=(x_i^\star-x_i)^T\Big(\nabla \phi_i(\alpha x_i^\star+(1-\alpha)x_i\big)-y_i\Big).
\end{align*}
Recalling
$x_i= \Pi_{\Omega_i}^{\phi_i}(y_i)=\mathop{\text{argmin}}_{x\in \Omega_i}\big\{-x^Ty_i+\phi_i(x)\big\}$, 
 $J(0)\leq J(\alpha)$, $\forall\alpha\in[0,1]$, because of the convexity of $\Omega_i$. This yields  $J'(\alpha)\big|_{0^+}\geq 0$, that is, 
\begin{align*}
J'(\alpha)|_{0^+}=(x_i^\star-x_i)^T(\nabla \phi_i(x_i)-y_i)\geq 0.
\end{align*}
On the other hand, 
\begin{align*}
(\bm \lambda-\bm \lambda^\star)^T(\bm \lambda-\bm \gamma)=([\bm\gamma]^+-\bm\lambda^\star)^T([\bm\gamma]^+-\bm\gamma)\leq0.
\end{align*}
Therefore, ${\beta}\leq 0$, which implies that the output $\bm z(t)$ is bounded.

Secondly, we show that $\bm s(t)$ is bounded. Actually, it follows from \eqref{V_1radially_unbdd} and the statement above that $\bm \gamma(t)$ is bounded. Thereby, we merely need to consider $\bm y$. Take another Lyapunov candidate function as 
\begin{align}\label{Lya_2}
V_2=\frac{1}{2}\|\bm y\|^2,
\end{align} 
which is radially unbounded in $\bm y$. Along the trajectories of Algorithm \ref{alg:1}, the derivative of $ V_2$ satisfies
\begin{align*}
\mathcal L_{\mathcal F}V_2=\big\{\zeta\in&\mathbb R:\zeta\in\sum_{i=1}^N y_i^T\big(-\partial f_i(x_i)-\partial g_i(x_i)^{\rm T}\lambda_i\\
 &-A_i^{\rm T}\mu_i+\nabla \phi_i(x_i)\big)-\|y_i\|^2\big\}.
\end{align*}
It is clear that $ \zeta\leq-\|\bm y\|^2+m\|\bm y\|=-2V_2+m\sqrt{2V_2}$ for a positive constant $m$, since $\bm x$, $\bm \lambda$ and $\bm \mu$ have been proved to be bounded. On this basis, 
it can be easily verified that $V_2$ is bounded, so is $\bm y$. Together, $\bm s(t)$ is bounded.

\noindent(ii) Set {$\mathcal R=\big\{(\bm z,\bm s):0\in \mathcal L_{\mathcal F}V _1\big\}$}. Clearly, by \eqref{optimal}, $\mathcal R\subseteq \big\{(\bm z,\bm s): \mathcal L(\bm x^\star,\bm \lambda,\bm \mu,\bm \omega^\star,\bm \nu^\star)=\mathcal L(\bm x,\bm \lambda^\star,\bm \mu^\star,\bm \omega,\bm \nu)\big\}$. Let $\mathcal M$ be the largest invariant subset of $\mathcal R$. By Lemma \ref{lem:inv_prin}, $(\bm z(t),\bm s(t))\rightarrow \mathcal M$ as $t\rightarrow \infty$.
Take any $(\tilde{\bm z},\tilde{\bm s})\in\mathcal M$. Let $\hat{\bm s}\in-F(\tilde{\bm z})+\nabla\Phi(\tilde{\bm z})$, and clearly $(\tilde{\bm z},\hat{\bm s})\in\mathcal M$ as well.  Similar to \eqref{Lya_1}, we take another Lyapunov function $\tilde V_1$  by replacing $(\bm z^\star,\bm s^\star)$ with $(\tilde{\bm z},\hat{\bm s})$. Based on similar arguments, $\tilde{\bm z}$ is Lyapunov stable, so is $\hat{\bm s}$. By Proposition 4.7 in \cite{haddad2011nonlinear}, there exists $(\bm z^\#,{\bm s}^\#)\in\mathcal M$ such that $({\bm z(t)},{\bm s(t)})\rightarrow ({\bm z}^\#,{\bm s}^\#)$ as $t\rightarrow \infty$, which yields that $\bm x(t)$ in Algorithm \ref{alg:1} converges to an optimal solution to problem \eqref{formulation}. \hfill$\square$

\begin{remark}
	It is worth mentioning that the {Bregman damping $\nabla\Phi(\bm z)$} in \eqref{compact} is fundamental to make the trajectory of variable $\bm s$ avoid going to infinity \cite{krichene2015accelerated,sun2020distributed}, or converging to an inexact optimal point \cite{gao2020continuous}. Clearly, $\bm z$ in the first ODE in \eqref{pro_alg} derives actually not from the variable $\bm z$ itself, but from the gradient of the quadratic function $\|\bm z\|^2/2$ instead. This is exactly the crucial point in designing the distributed MD-based dynamics \eqref{compact}. Correspondingly, the properties in conjugate spaces, referring to \eqref{fenchel}, play an important role in the analysis.
\end{remark}

%
%

%
For convenience, we define
\begin{align*}
\hat{\bm x}\triangleq& \frac{1}{t}\int_{0}^{t}\!\bm x(\tau)d\tau,\;\hat{\bm \lambda}\triangleq \frac{1}{t}\!\int_{0}^{t}\bm \lambda(\tau)d\tau,\; \hat{\bm \mu}\triangleq \frac{1}{t}\!\int_{0}^{t}\bm \mu(\tau)d\tau,\\
\hat{\bm \omega}\triangleq& \frac{1}{t}\int_{0}^{t}\bm \omega(\tau)d\tau,\;\hat{\bm \nu}\triangleq \frac{1}{t}\int_{0}^{t}\bm \nu(\tau)d\tau.
\end{align*}
Then we describe the convergence rate of Algorithm \ref{alg:1}.
\begin{theorem}
Under  Assumption \ref{ass_1}, \eqref{compact} converges with a rate of $\mathcal O(1/t)$, \textit{i.e.},
\begin{align*}
0\leq\mathcal L(\hat{\bm x},\bm \lambda^\star,\bm \mu^\star,\hat{\bm \omega},\hat{\bm \nu})-\mathcal L(\bm x^\star,\hat{\bm \lambda},\hat{\bm \mu},\bm \omega^\star,\bm \nu^\star)\leq \frac{1}{t}V_1(\bm s(0)).
\end{align*}
\end{theorem}
\textbf{Proof}. 
It follows from \eqref{deriv_V1}-\eqref{deriv_V11} that 
\begin{align*}
\frac{d}{dt}V_1\leq\mathcal L(\bm x^\star,\bm \lambda,\bm \mu,\bm \omega^\star,\bm \nu^\star)-\mathcal L(\bm x,\bm \lambda^\star,\bm \mu^\star,\bm \omega,\bm \nu)\leq 0.
\end{align*}
By integrating both sides over the time interval $[0,t]$, 
\begin{align}\label{rate}
-V_1(\bm s(0))\leq V_1(\bm s(t))-V_1(\bm s(&0))\nonumber\\
\leq\int_{0}^t\Big(\mathcal L(\bm x^\star,\bm \lambda(\tau),\bm \mu(\tau),\bm \omega^\star,\bm \nu^\star)&\\
-\mathcal L(\bm x(\tau),\bm \lambda^\star,\bm \mu^\star,&\bm \omega(\tau),\bm \nu(\tau))\Big)d\tau\leq 0.\nonumber
\end{align}
With applying Jensen's inequality to the convex-concave  Lagrangian function  $\mathcal L$,
\begin{align*}
\mathcal L(\bm x^\star\!,\bm \lambda(t),\bm \mu(t),\bm \omega^\star\!,\bm \nu^\star)\!&\geq\! \frac{1}{t}\!\int_{0}^t\!\!\mathcal L(\bm x^\star\!,\bm \lambda(\!\tau\!),\bm \mu(\!\tau\!),\bm \omega^\star\!,\bm \nu^\star)d\tau,\\\mathcal L(\bm x(t),\!\bm \lambda^\star\!,\bm \mu^\star\!,\bm \omega(t),\!\bm \nu(t)\!)\!&\leq\! \frac{1}{t}\!\int_{0}^t\!\!\mathcal L(\bm x(\!\tau\!),\!\bm \lambda^\star\!,\bm \mu^\star\!,\bm \omega(\!\tau\!),\!\bm \nu(\!\tau\!)\!)d\tau.
\end{align*}
By substituting the above inequalities into \eqref{rate}, the conclusion follows. \hfill$\square$

\section{Numerical examples}\label{Sec:num}

\begin{figure}
	\centering
	\subfigure[Trajectories of $x_i$]{
		\includegraphics[width=4.1cm]{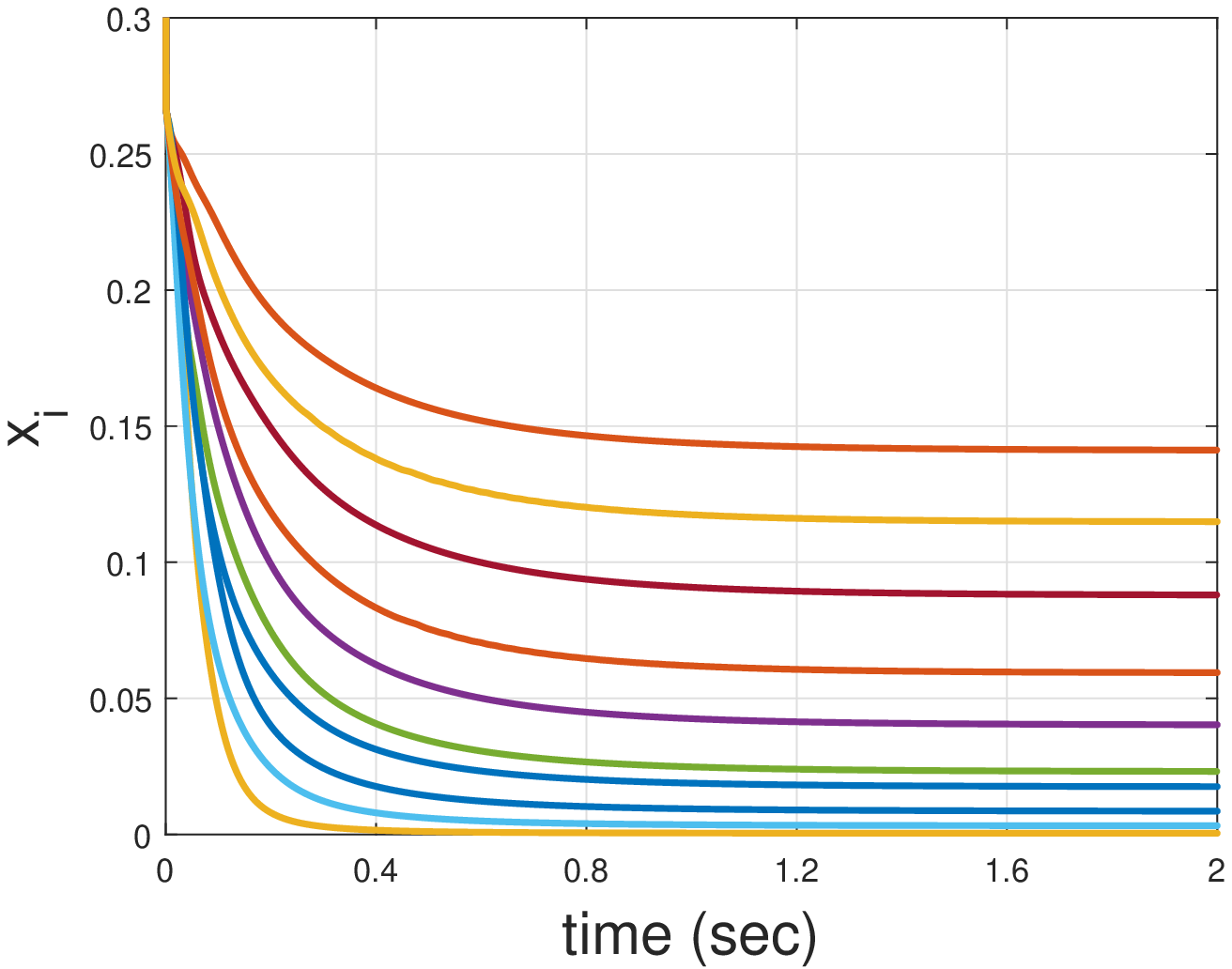}
	}
	\,\hspace{-0.4cm}
	\subfigure[Trajectories of $y_i$]{
		\includegraphics[width=4.1cm]{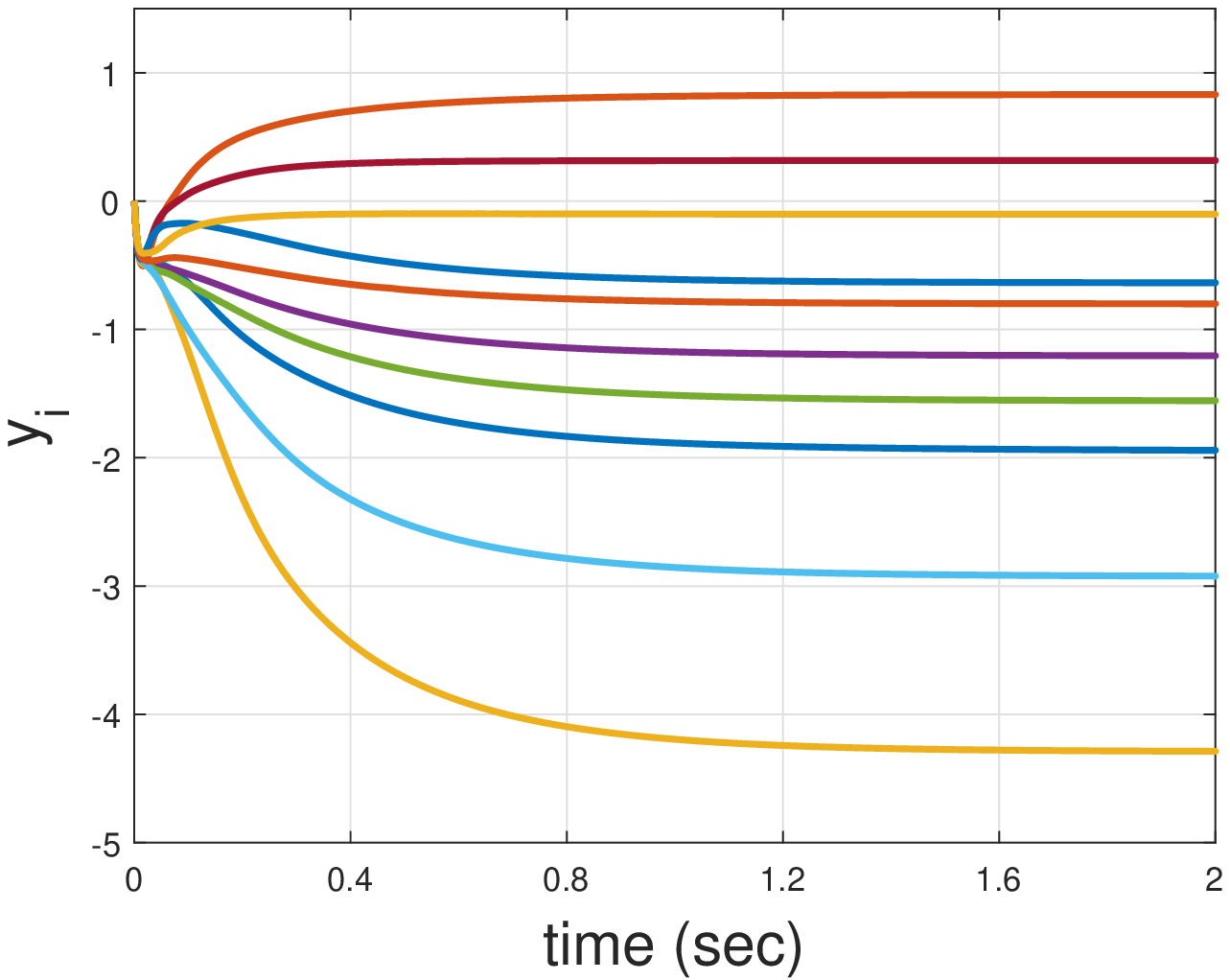}
	}
	\caption{Trajectories of all agents' variables.}
	\label{fig6}
\end{figure}

In this section, we examine the correctness and effectiveness of Algorithm 1 on the classical simplex-constrained problems (see, e.g., \cite{gao2020continuous,yuan2018optimal}), where  the local constraint set is an  $ n $-simplex, e.g., $$ \Omega_{i}=\{x_{i} \in \mathbb{R}^{n}_{+}: \sum_{k=1}^{n} x_{i,k}=1\},\, \forall i \in \mathcal V .$$
First, we consider  the following nonsmooth optimization problem with $ N = 10 $ and  $ n = 4 $,  
\begin{equation}\label{f3}
	\begin{aligned}
		& \min _{\boldsymbol{x}\in\boldsymbol{\Omega}}\, \sum_{i=1}^{N} \left\| W_{i}x_{i}-d_{i}\right\|^{2}+c_{i}\left\|x_{i}\right\|_{1} \\
		\text { s.t. } & 	\sum_{i=1}^{N} g_{i}(x_{i})\leq 0, \quad	
		\sum_{i=1}^{N} A_{i} x_{i}-\sum_{i=1}^{N} b_{i}=0_{2},
	\end{aligned}
\end{equation}
where $ W_{i}   $ is a positive semi-definite matrix,  $ d_{i} \in  \mathbb{R}^{4} $,  and $c_{i}>0$.  The coupled inequality constraint is
$$
g_{i}(x_{i})=\left\|x_{i}\right\|^{2}+c_{i}\left\|x_{i}\right\|_{1}- \frac{25}{2n+i^{2}}\; ,
$$ 
while $ A_{i} \in \mathbb{R}^{2\times4} $  and $ b_{i} \in \mathbb{R}^{2}$ are random matrices ensuring  the Slater’s constraint condition. Here, $ W_{i} $, $ d_{i} $, $ g_{i} $, $ A_{i} $ and $ b_{i} $ are  private to agent $ i $, and all agents   communicate through an undirected cycle
network $ \mathcal{G} $:
$$
1\rightleftarrows 2\rightleftarrows\cdots\rightleftarrows10\rightleftarrows1.
$$
To implement the MD method, we employ  the  negative
entropy function $ \phi_i(x_i)\!=\!\sum_{k=1}^nx_{i,k}\!\log(x_{i,k}) $ as the  generating function on $ \Omega_{i} $ in Algorithm \ref{alg:1}.
In Fig. \ref{fig6}, we show the  trajectories of one dimension of each  $ x_{i} $ and $y_i$, respectively.
Clearly,  the trajectories of both $x_i$ and $y_i$ in MDBD are bounded, while the boundedness of $y_i$ may not be guaranteed in \cite{krichene2015accelerated,sun2020distributed}.

Next, we show the effectiveness of MDBD by comparisons. As is investigated in \cite{ben2001ordered,krichene2015accelerated}, when the generating function satisfies $ \phi_i(x_i)\!=\!\sum_{k=1}^nx_{i,k}\!\log(x_{i,k}) $ on the unit simplex, $ \Pi_{\Omega_{i}}^{\phi_{i}}(y_{i}) $ can be explicitly expressed as  \eqref{qq}.
%
In this circumstance, the MD-based method works better than projection-based algorithms, since it can  be regarded as projection-free and effectively saves the time for projection operation,  especially with high-dimensional variables. 

\begin{figure}[htbp]
	\centering
	\subfigure[$ n=4 $]{
		\includegraphics[width=4.1cm]{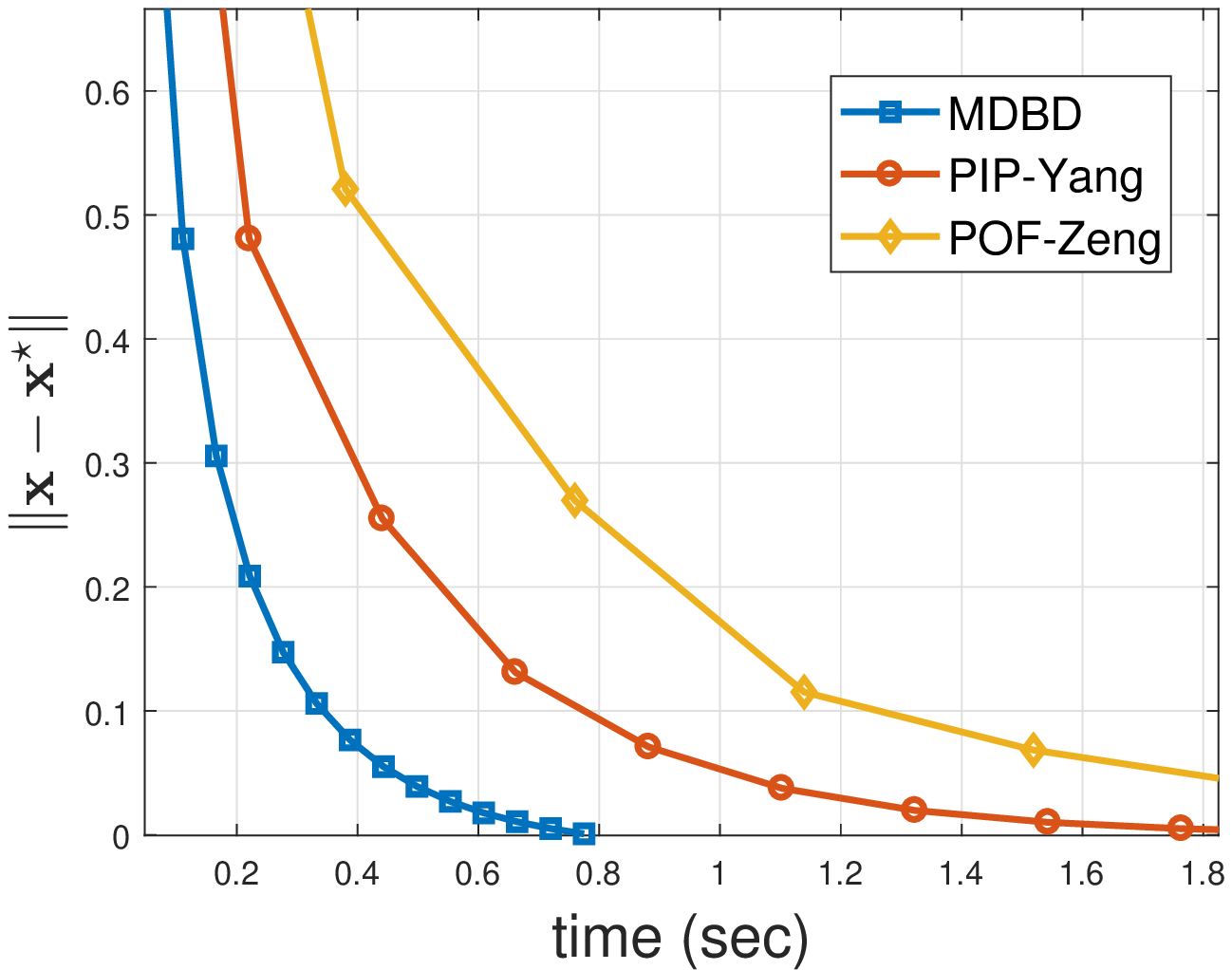}
	}
	\,\hspace{-0.4cm}
	\subfigure[$ n=16 $]{
		\includegraphics[width=4.1cm]{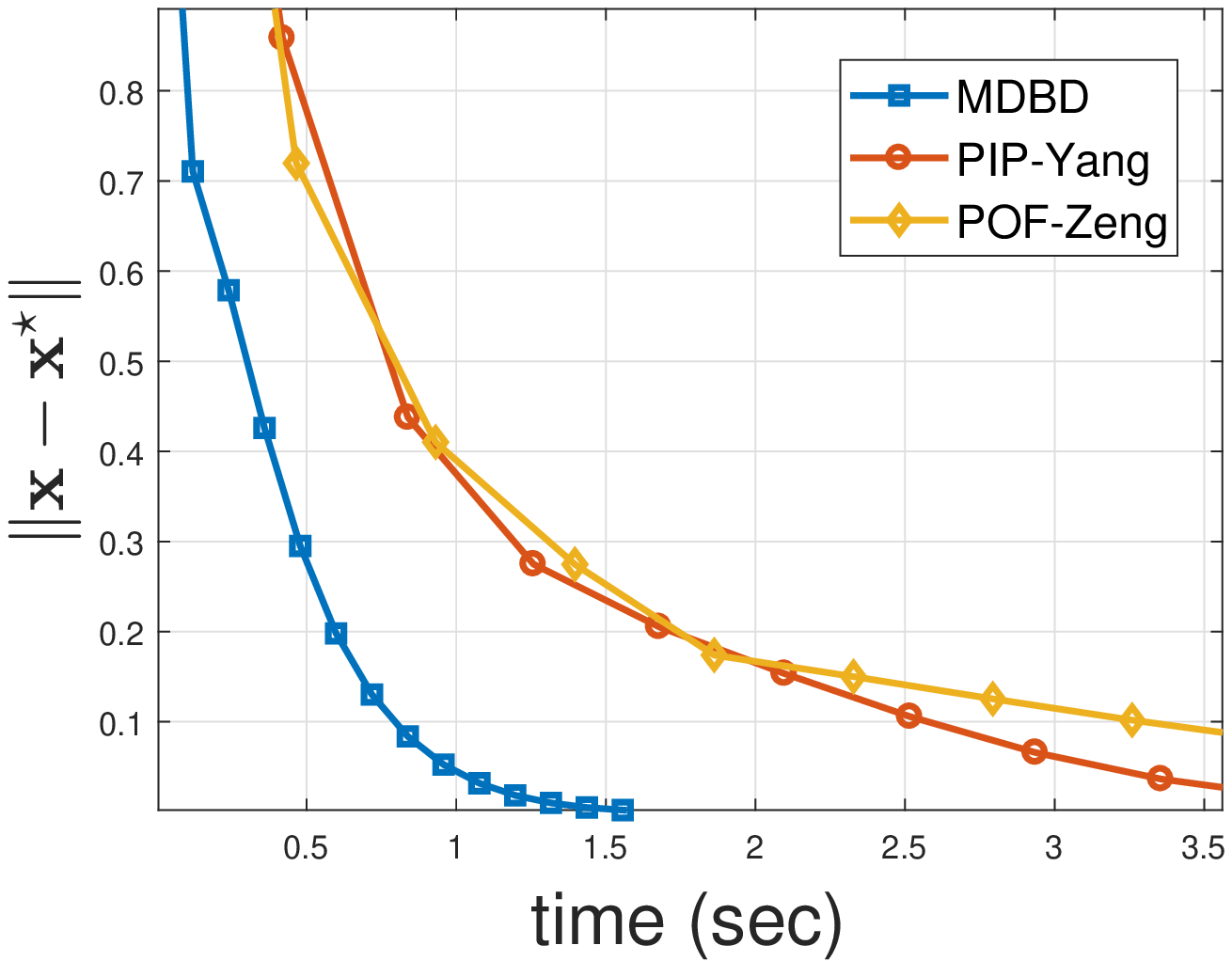}
	}
	\,
	\subfigure[$ n=64 $]{
		\includegraphics[width=4.1cm]{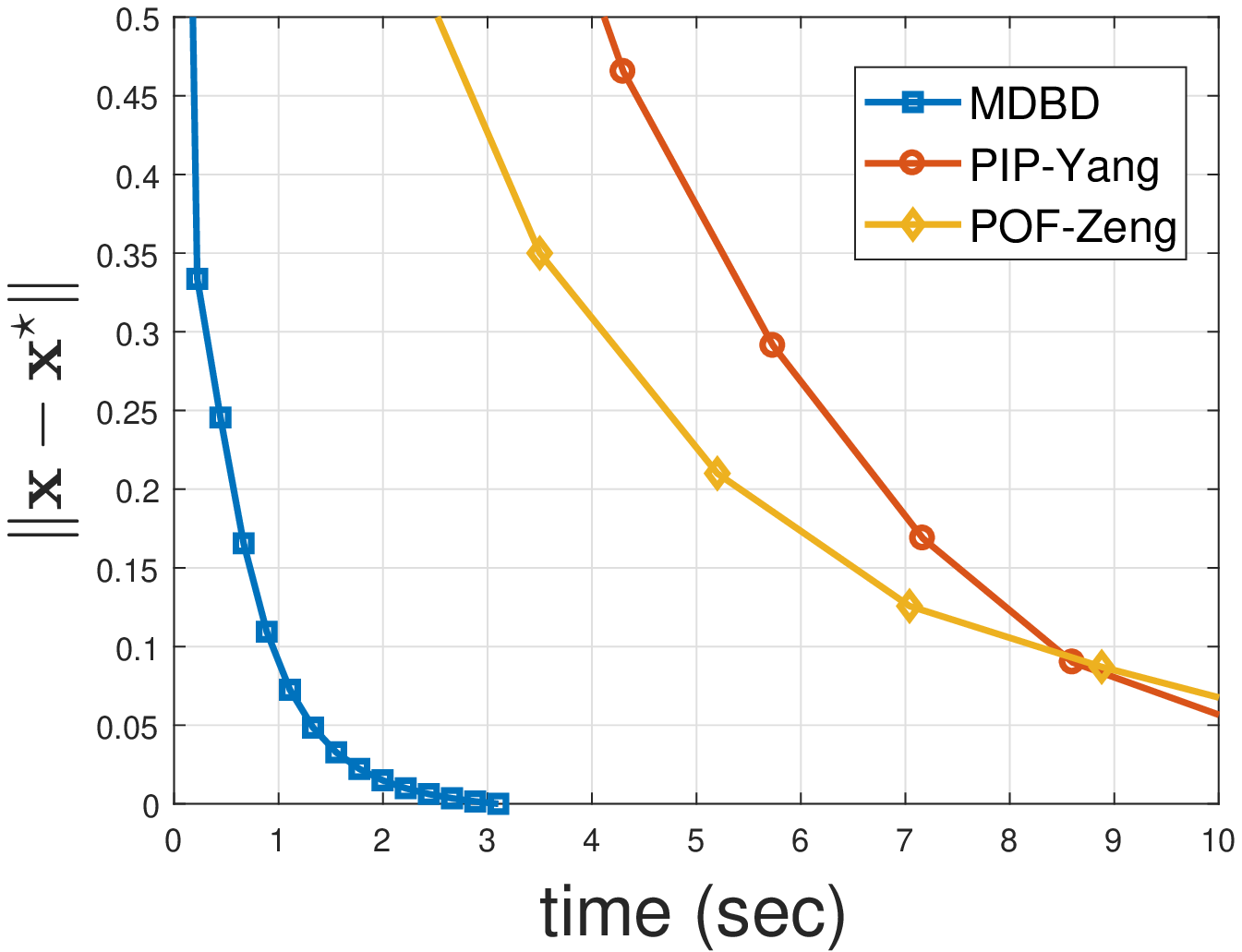}
	}
	\,\hspace{-0.4cm}
	\subfigure[$ n=256 $]{
		\includegraphics[width=4.1cm]{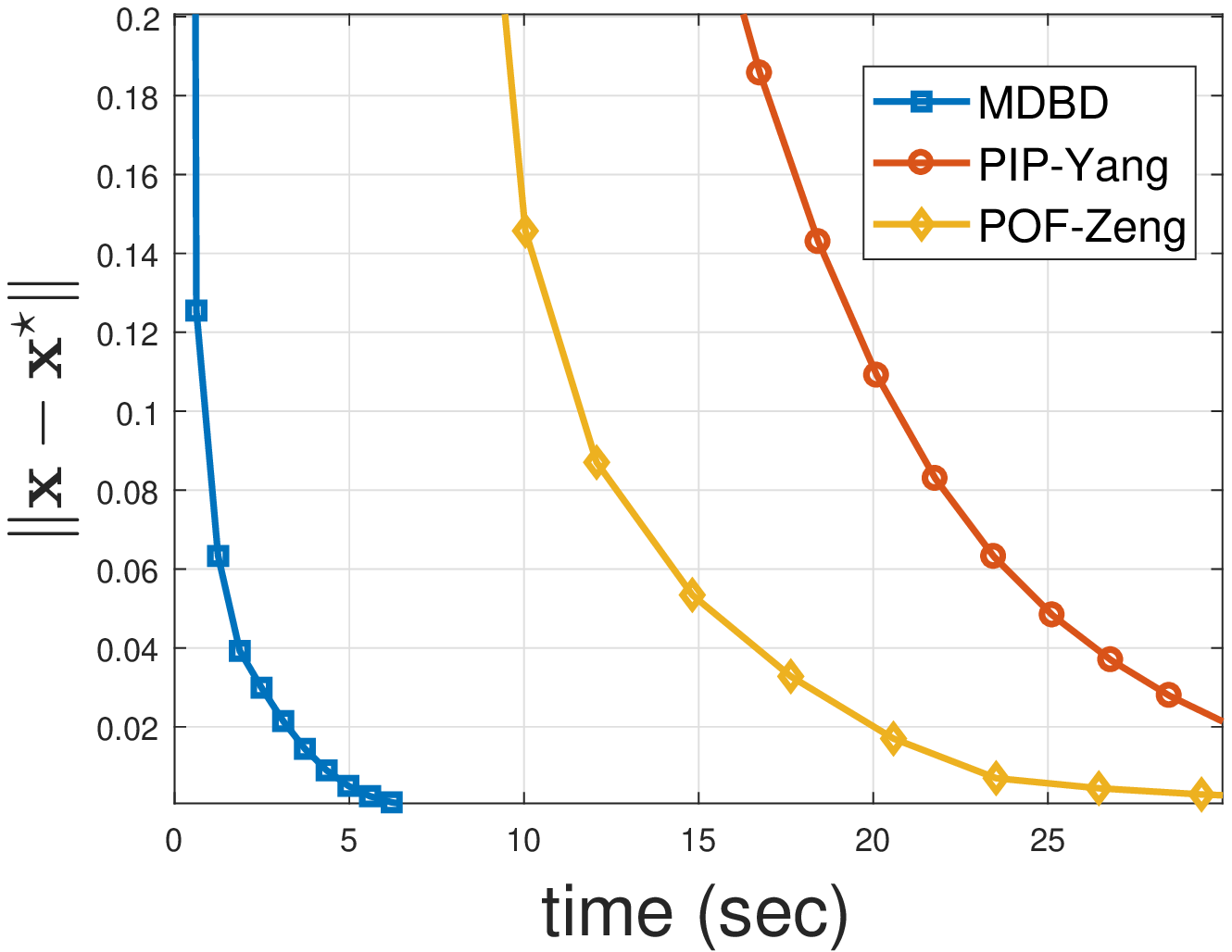}
	}
	\caption{Optimal errors in different dimensions: $ n =
		4, 16, 64, 256 $.}
	\label{fig3}
\end{figure}

To this end,  we investigate different dimensions of
decision variable $x_i$ and compare  MDBD with two distributed continuous-time projection-based algorithms --- the  proportional-integral protocol (PIP-Yang) in \cite{yang2016multi} and the  projected output feedback (POF-Zeng) in \cite{zeng2018distributed}, still for the cost
functions and the coupled constraints given in
\eqref{f3}. 
\begin{table*}
	\centering
	\small
	\caption{Real running time (sec) in different dimensions} 
	\setlength\tabcolsep{12pt}
	\renewcommand\arraystretch{1.2}
	\begin{tabular}{c|c|c|c|c|c|c|c}
		\hline
		\hline
		\specialrule{0em}{0.2pt}{0.5pt}
		 & 	 $ n=4 $&  $ n=64 $ & $ n=256 $& $n=1024$ &$ n=4096 $& $ n=10^{5} $& $ n=10^{6}$\\ \hline	
		MDBD &0.47 &2.42 &6.76 &12.98  &27.99 &146.62  & 466.60\\ \hline
		PIP-YANG&2.51&19.63 &48.51 & 195.67& 892.74 &$ >3000 $ &$ >5000 $\\  \hline
		POF-ZENG&3.92 &21.78 &39.73 &207.03& 1136.85 & $ >3000 $  &$ >5000 $ \\\hline
		\hline
	\end{tabular}
	\label{tab1}
\end{table*}
In Fig. \ref{fig3}, the $ x $-axis is for the \textit{real running time} of the GPU, while the $ y $-axis is for the optimal error $ \| \bm x-\bm x^\star\| $. As the dimension increases, the \textit{real running time} of the two projection-based dynamics is obviously longer than that of MDBD, because obtaining \eqref{qq}    is much faster than calculating a projection on high-dimensional constraint sets via solving a general quadratic optimization problem.

Furthermore, Table \ref{tab1} lists the \textit{real running time} for
three algorithms with different dimensions of decision variables. As the dimension increases, finding the projection points in large-scale circumstances becomes more and more difficult. But remarkably, MDBD still maintains good performance, due to the advantage of MD. 

\section{CONCLUSIONS}\label{Sec:con}
We investigated distributed nonsmooth optimization with both local set constraints and coupled constraints. Based on the mirror descent method, we proposed a continuous-time algorithm with introducing the Bregman damping  to guarantee the algorithm's boundedness and accuracy. Furthermore, we utilized nonsmooth techniques, conjugate functions, and the Lyapunov stability to prove the convergence. Finally, we implemented comparative experiments to illustrate the effectiveness of our algorithm.



%
%
%
%
\bibliographystyle{IEEEtran}
\bibliography{reference}
\end{document}